\documentclass[12pt]{article} 

\usepackage{latexsym}\usepackage{amssymb,latexsym} \usepackage{subfig}
\usepackage[mathscr]{eucal} \usepackage{amsthm} \usepackage{amsmath}
\usepackage{amsfonts} \usepackage{layout}
\usepackage{graphicx}
\usepackage{bold-extra}
\usepackage{float} \restylefloat{figure}
\usepackage{multirow}

\begin{document}

\title{\bf Confidence Bands for the Logistic and Probit Regression
  Models Over Intervals}\author{Lucy Kerns\\ Department of Mathematics
  and Statistics\\ Youngstown State University\\ Youngstown OH 44555
  \\xlu@ysu.edu\date{} }

\maketitle

\noindent \textbf{Summary}

\noindent This article presents methods for the construction of
two-sided and one-sided simultaneous hyperbolic bands for the logistic
and probit regression models when the predictor variable is restricted
to a given interval. The bands are constructed based on the asymptotic
properties of the maximum likelihood estimators. Past articles have
considered building two-sided asymptotic confidence bands for the
logistic model, such as Piegorsch and Casella (1988). However, the
confidence bands given by Piegorsch and Casella are conservative under
a single interval restriction, and it is shown in this article that
their bands can be sharpened using the methods proposed
here. Furthermore, no method has yet appeared in the literature for
constructing one-sided confidence bands for the logistic model, and no
work has been done for building confidence bands for the probit model,
over a limited range of the predictor variable.  This article provides
methods for computing critical points in these areas.  \bigskip


\newpage \setcounter{page}{1}

\bigskip

\section {Introduction}
Logistic and probit regression models are widely used for modeling
dichotomous outcomes, and have been increasingly applied in medical
research, public health research, environmental science, and many
other fields, such as human behavior modeling (Chou, Lu, and Mao,
2002), environmental modeling (Pradhan and Lee, 2010), and biomedical
research (Austin and Steyerberg, 2014).  The logistic and probit
regression models are statistical methods that allow one to estimate
the response probability for a dichotomous response, that is, a
response which is binary, taking values 1 (success, normal, positive,
etc.) and 0 (failure, abnormal, negative, etc.).  In this article, we
consider the case where the binary response variable $Y$ is determined
by a predictor variable $x$, and the response probability is:
\[P(Y = 1) = p(x)=
 \begin{cases}
      1/[1 + exp(-\boldsymbol{c}'\boldsymbol{\beta})], & \text{logistic model} \\
      \Phi(\boldsymbol{c}'\boldsymbol{\beta}),  & \text{probit model}
 \end{cases}
\]
where $\boldsymbol{c} = (1 \quad x)'$,
$\boldsymbol{\beta}=(\beta_0 \quad \beta_1)'$, and $\Phi$ is the
cumulative distribution function (cdf) of the standard normal
distribution.

The construction of confidence bands on $p(x)$ is often of interest,
and we will use the following link functions to transform the problem
of constructing confidence bands for the response probability $p(x)$
to the problem of constructing confidence bands for the linear
predictor $\boldsymbol{c}'\boldsymbol{\beta}$ on which the bands are
defined.
\[
\text{logit}(p(x))= \text{log}_e\left[\frac{p(x)}{1-p(x)}\right]= \boldsymbol{c}'\boldsymbol{\beta}, \,\, \text{logistic regression}
\]
\[
\text{probit}(p(x))= \Phi^{-1}(p(x))= \boldsymbol{c}'\boldsymbol{\beta}, \,\, \text{probit regression}
\]

There is a wealth of literature on building exact or conservative
two-sided bands for the linear predictor
$\boldsymbol{c}'\boldsymbol{\beta}$ with one or more than one
predictor variable. This includes the work of Scheff\'{e} (1953), and
Working $\&$ Hotelling (1929), among others. Many authors have
improved earlier work by restricting the predictor variables to given
intervals, including Wynn and Bloomfield (1971), Casella and
Strawderman (1980), and Uusipaikka (1983). Liu \emph{et al.} (2005)
developed a simulation-based method for constructing two-sided
confidence bands over a restricted region for the multiple regression
model. Less work has been done on the construction of one-sided (lower
or upper) confidence bands for linear regression. With no
restriction on the predictor variable, Hochberg and Quade (1975)
developed a method for constructing one-sided bands in the multiple
regression setting. When the predictor variable is constrained to a
pre-specified interval, Bohrer and Francis (1972) gave exact one-sided
hyperbolic confidence bands for the simple linear case. Pan \emph{et
al.} (2003) extended Uusipaikka's work (1983) to the computation of
one-sided bands under a single interval restriction.  Liu \emph{et al.}'s
article (2008) summarized several existing methods and also provided
new methods for the construction of two-sided and one-sided confidence
bands with a restricted predictor variable.

Much less work has been done on constructing confidence bands for the
logistic regression model. Brand, Pinnock, and Jackson (1973)
constructed confidence bands for both $p(x)$ and the inverse of $p(x)$
when there is one predictor variable and no constraints exist on the
predictor variable. Hauck (1983) extended their work to more than one
predictor variable yet still with no constraints. In the case of
restricted predictor variables, Piegorsch and Casella (1988) extended
the work of Casella and Strawderman (1980) from the multiple
regression model to the logistic case. The method of Casella and
Strawderman produced asymptotic two-sided bands over a particular
constrained region of predictor variables. However, if the restricted
region on predictor variables is rectangular, the resulting confidence
bands are rather conservative even with large samples.

Wei Liu's book (2010, Chapter 8) presented confidence bands for the logistic model
with more than one explanatory variable. The method utilized
simulation-based confidence bands (Liu \emph{et al.} (2005)) for the
linear predictor $\boldsymbol{c}'\boldsymbol{\beta}$ in the multiple
linear model, and the desired bands for the logistic model were then
obtained via the logit link function. The new contribution of this paper is that we provide explicit expressions for determining the critical values in the logistic setting. Wei Liu's method is very broad, but relies on simulation. Our method is more focused but admits tractable forms.

Furthermore, no work has yet appeared for building one-sided
confidence bands for the logistic regression model when the predictor
variable is bounded on a given interval, and no methods have been
developed for constructing both two-sided and one-sided bands for the
probit regression model with a restricted predictor variable. In this
paper, we center our attention on building asymptotic two-sided and
one-sided hyberbolic confidence bands for the logistic and probit
models over a limited range of the predictor variable, say, between
$(a,\, b)$, where $a$ and $b$ are constants.

For the logistic and probit models, we denote the ML estimator of
$\boldsymbol{\beta}$ by $\hat{\boldsymbol{\beta}}$, and under certain
regularity conditions (Kendall and Stuart (1979)),
$\hat{\boldsymbol{\beta}}$ follows asymptotically as:
\[
\hat{\boldsymbol{\beta}} \sim ASN_2(\boldsymbol{\beta},\boldsymbol{F^{-1}}),
\]
where $\boldsymbol{F}$ is the Fisher information matrix. It is well
known that the logistic model has an explicit formula for the Fisher
information matrix, while the probit model does not, but the
information matrix can be obtained through numerical methods.

For $x \in (a,b)$, where $a$ and $b$ are given constants, a
$100(1-\alpha)\%$ two-sided hyperbolic band for the linear function
$\boldsymbol{c}'\boldsymbol{\beta}$ has the form
\begin{equation*}
\boldsymbol{c'\beta} \in \boldsymbol{c'\hat{\beta}} \pm w(\boldsymbol{c'F^{-1}c})^{1/2},
\end{equation*}
where $\boldsymbol{c} = (1 \quad x)'$ and $w$ satisfies
\begin{equation}\label{eq2.1}
P[\boldsymbol{c'\beta} \in \boldsymbol{c'\hat{\beta}} \pm w(\boldsymbol{c'F^{-1}c})^{1/2},\,\,\text{for all} \,\, x \in (a,b)] = 1 - \alpha.
\end{equation}

An upper $100(1-\alpha)\%$ one-sided hyperbolic band for the linear
function $\boldsymbol{c}'\boldsymbol{\beta}$ has the form
\begin{equation*}
\boldsymbol{c'\beta} \leq \boldsymbol{c'\hat{\beta}} + w_u(\boldsymbol{c'F^{-1}c})^{1/2},
\end{equation*}
where $w_u$ satisfies
\begin{equation}\label{eq2.2}
P[\boldsymbol{c'\beta} \leq \boldsymbol{c'\hat{\beta}} + w_u(\boldsymbol{c'F^{-1}c})^{1/2},\,\, \text{for all}\,\, x \in (a,b)] = 1 - \alpha.
\end{equation}
A $100(1-\alpha)\%$ lower one-sided hyperbolic band for
$\boldsymbol{c}'\boldsymbol{\beta}$ can be defined similarly.

Since both logit and probit functions are monotonically increasing
functions, $100(1-\alpha)\%$ two-sided confidence bands for $p(x)$ in
the logistic regression model are given by
\begin{equation*}
 \{1 + \text{exp}[- \boldsymbol{c'\hat{\beta}} + w(\boldsymbol{c'F^{-1}c})^{1/2}]\} ^{-1} \leq p(x) \leq
 \{1 + \text{exp}[- \boldsymbol{c'\hat{\beta}} - w(\boldsymbol{c'F^{-1}c})^{1/2}]\} ^{-1},
\end{equation*}
and $100(1-\alpha)\%$ two-sided confidence bands for $p(x)$ in the
probit regression model are given by
\begin{equation*}
 \Phi(\boldsymbol{c'\hat{\beta}} - w(\boldsymbol{c'F^{-1}c})^{1/2}) \leq p(x) \leq
 \Phi(\boldsymbol{c'\hat{\beta}} + w(\boldsymbol{c'F^{-1}c})^{1/2}).
\end{equation*}
Similarly, $100(1-\alpha)\%$  upper confidence bands for $p(x)$ are given by
\[p(x) \leq
\begin{cases}
  \{1 + \text{exp}[- \boldsymbol{c'\hat{\beta}} - w_u(\boldsymbol{c'F^{-1}c})^{1/2}]\}^{-1}, &\text{logistic model}\\
  \Phi(\boldsymbol{c'\hat{\beta}} + w_u(\boldsymbol{c'F^{-1}c})^{1/2}), &\text{probit model}
\end{cases}
\]

In what follows, we consider the logistic model and focus on
Equations~(\ref{eq2.1}) and (\ref{eq2.2}). In particular, we
propose methods to find the critical values $w$ in
Equation~(\ref{eq2.1}) for the two-sided case and $w_u$ in
Equation~(\ref{eq2.2}) for the one-sided case under the logistic
model. The derivation of the critical values under the probit model is
the same except $\boldsymbol{F}^{-1}$ is replaced by
$\boldsymbol{\Omega}^{-1}$.

Both the logistic and probit models are special cases of the
generalized linear model (GLM), and the methodology proposed here can
also be applied to other forms of GLM (the complementary-log-log
model, for example), that can be transformed into the standard
regression model via a link function. In each case, a set of
simulations are required to confirm the validity of the method in
small samples, thus in the interest of brevity only logistic and
probit models are considered here. We devote most of our attention to
the logistic model to illustrate the methodology.

Liu \emph{et al.}'s article (2008) presents several methods for the
construction of two-sided and one-sided confidence bands for the
simple linear model with a restricted predictor variable. Since we can
transform the logistic and probit regression setting to the simple
linear regression setting via the link functions (logit and probit),
we extend their work to the logistic and probit regression models and
develop methods to derive asymptotic confidence bands for large
samples. This paper is organized as follows. We give the general
setting of the problem in Section 1. Two-sided confidence bands are
discussed in Section 2, and Section 3 presents the results in the
one-sided case. In Section 4, a Monte Carlo simulation is run to
investigate how well the asymptotic approximation holds for small
sample sizes. An example is given in Section 5 to illustrate the
proposed methods.  \bigskip

\section{Two-Sided Bands}
Our construction of confidence bands are based on the methods given by
Liu \emph{et al.} (2008). There are undoubtedly some similarities between the methods presented here and in their review. The main difference between their method and ours
lies at the beginning of the derivation. The formulation in their review relied on a standard bivariate \emph{t} random variable (See
Equations~(13) and (14) in their review), while the derivation given
here involves a standard \emph{z} bivariate random variable, as will be soon shown.  As a result of this reformulation, our critical values are based on a chi-square random variable, instead of an $F$ variate.  In the interest of completeness, we include most of the computational
details regarding the critical values $w$ and $w_u$.

The Fisher information matrix $\boldsymbol{F^{-1}}$ is a positive
semi-definite symmetric matrix, so there exists a positive
semi-definite matrix $\boldsymbol{B}$ such that
$\boldsymbol{F^{-1}} = \boldsymbol{B}^2$. Then
\[
\boldsymbol{z} = \boldsymbol{B}^{-1}(\boldsymbol{\beta} - \hat{\boldsymbol{\beta}}) \sim N_2(\boldsymbol{0},\boldsymbol{I}).
\]
Define the polar coordinates of $\boldsymbol{z} = (z_1,z_2)'$, $(R_z,Q_z)$, by
\[
z_1 = R_z\,\mathrm{cos}Q_z,\,\,z_2 = R_z\,\mathrm{sin}Q_z,\,\,\text{for}\,\, R_z \geq 0 \,\,\text{and}\,\, Q_z \in [0,2\pi].
\]
It is well known that $R_z^2$ has the $\chi_2^2$ distribution and is statistically independent of $Q_z$, which has a uniform $[0,2\pi]$ distribution.
Write Equation~(\ref{eq2.1}) as
\begin{align}\label{eq2.3}
 &P[\boldsymbol{c'\beta} \in \boldsymbol{c'\hat{\beta}} \pm w(\boldsymbol{c'F^{-1}c})^{1/2}, \,\,\text{for all}\,\, x \in (a,b)] \notag\\
 &=P\left[\sup\limits_{x \in (a,b)}
     \frac{|\boldsymbol{c'\beta} - \boldsymbol{c'\hat{\beta}}|}{(\boldsymbol{c'F^{-1}c})^{1/2}} < w \right]\notag\\
 &=P\left[\sup\limits_{x \in (a,b)}
     \frac{|\boldsymbol{(Bc)'z}|}{\lVert\boldsymbol{Bc}\rVert} < w \right]
\end{align}
Notice that for a given vector $\boldsymbol{u} \in R^2$ and a number
$w > 0$,
$\{\boldsymbol{z}:\boldsymbol{u}'\boldsymbol{z}/\lVert\boldsymbol{u}\rVert
= w\}$
represents a straight line that is perpendicular to the vector
$\boldsymbol{u}$ and is on the same side of the origin as the vector
$\boldsymbol{u}$. The perpendicular distance from the origin to this
line is $w$. Therefore, the set defined by
\begin{equation*}
   \{\boldsymbol{z}:|\boldsymbol{u}'\boldsymbol{z}|/\lVert\boldsymbol{u}\rVert < w \} \subset R^2
\end{equation*}
consists of all points that are sandwiched between two parallel
straight lines
$\boldsymbol{u}'\boldsymbol{z}/\lVert\boldsymbol{u}\rVert = w$ and
$\boldsymbol{u}'\boldsymbol{z}/\lVert\boldsymbol{u}\rVert =
-w$.
Therefore, letting $\boldsymbol{Bc} = \boldsymbol{u}$,
Equation~(\ref{eq2.3}) can be further expressed as:
\begin{equation}\label{eq2.4}
 P\left[\sup\limits_{x \in (a,b)}
     \frac{|\boldsymbol{(u)'z}|}{\lVert\boldsymbol{u}\rVert} < w \right]
 =P\{\boldsymbol{z} \in R_2\},
\end{equation}
where $R_2 = \cap_{x \in (a,b)} R_2(x),$ and
$R_2(x) =
\{\boldsymbol{z}:|\boldsymbol{u}'\boldsymbol{z}|/\lVert\boldsymbol{u}\rVert
< w\} $.
$R_2$ is depicted in Figure~1(a). We can rotate the region $R_2$
around the origin so that the angle $\phi$ between the two vectors
$\boldsymbol{u_a}=\boldsymbol{B}(1 \quad a)$ and
$\boldsymbol{u_b}=\boldsymbol{B}(1 \quad b)$ is equally divided by the
$z_1$-axis. The new region, denoted by $R_2^*$, is depicted in Figure
1(b). Because of the rotation invariance of the normal distribution,
we have $P\{\boldsymbol{z} \in R_2\} = P\{\boldsymbol{z} \in R_2^*\}$,
and $R_2^*$ can be expressed as
\[
  R_2^* = \{\boldsymbol{z}:|\boldsymbol{u}'\boldsymbol{z}|/\lVert\boldsymbol{u}\rVert < w, \text{for all}\,\,
  \boldsymbol{u} \in E(\phi)\},
\]
where $E(\phi)$ is a cone depicted as the shaded area in Figure 1(b),
and
$E(\phi) = \{\boldsymbol{u}: u_2 > \lVert\boldsymbol{u}\rVert
\mathrm{cos}(\phi/2)\}$.
Therefore, the simultaneous confidence level is equal to
\begin{equation}\label{eq2.5}
  P\left[\sup\limits_{\boldsymbol{u} \in E(\phi)}
     \frac{|\boldsymbol{u'z}|}{\lVert\boldsymbol{u}\rVert} < w \right]
 =P\{\boldsymbol{z} \in R_2^*\}.
\end{equation}
\begin{figure}[H]
\centering
\subfloat[Region $R_2$]{\includegraphics[width=3.5in,height=2.9in]{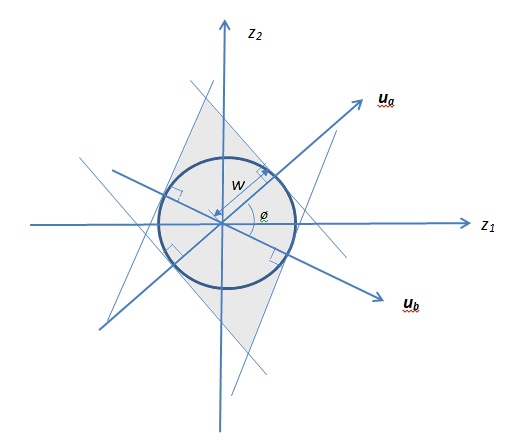}}
\subfloat[Region $R_2^*$]{\includegraphics[width=3.8in,height=2.9in]{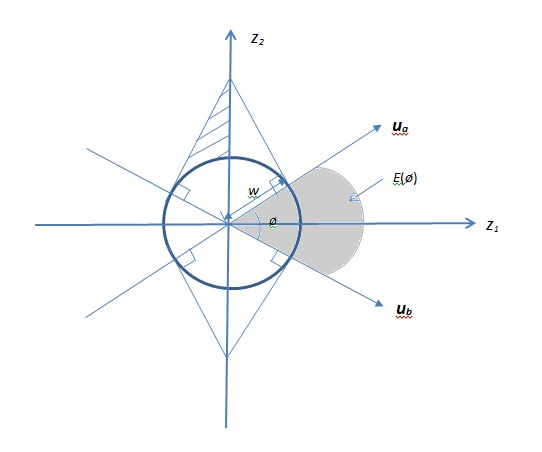}}
 \caption{The Regions $R_2$ and $R_2^*$}
\end{figure}

\subsection{Method 1: Finding Expression for the Supremum}
This method is based on finding an exact expression for the
supremum. Notice that
\[
\frac{|\boldsymbol{u'z}|}{\lVert\boldsymbol{u}\rVert}
=  \lVert\boldsymbol{z}\rVert\frac{|\boldsymbol{u'z}|}{\lVert\boldsymbol{u}\rVert\lVert\boldsymbol{z}\rVert}
= \lVert\boldsymbol{z}\rVert \lvert \mathrm{cos}(\theta_{uz})\rvert,
\]
where $\theta_{uz}$ is the angle between $\boldsymbol{u}$ and
$\boldsymbol{z}$. Denote
$\Psi_1=[-\phi/2,\phi/2] \cup [\pi-\phi/2,\pi+\phi/2]$,
$\Psi_2=[\phi/2,\pi/2] \cup [\pi+\phi/2,3\pi/2]$, and
$\Psi_3=[\pi/2,\pi-\phi/2 ] \cup [3\pi/2,-\phi/2]$.  Since
$\mathrm{cos}(\theta)$ is monotonically decreasing on $[0, \pi]$, we
have the following results:
\[
\sup_{\boldsymbol{u} \in E(\phi)}\lvert \mathrm{cos}(\theta_{uz})\rvert =
 \begin{cases}
      1                                    &\text{if}\,\, \theta_z \in \Psi_1, \\
      \lvert \mathrm{cos}(\theta_z - \phi/2)\rvert  & \text{if}\,\, \theta_z \in \Psi_2,\\
      \lvert \mathrm{cos}(\theta_z + \phi/2)\rvert  & \text{if}\,\, \theta_z \in \Psi_3.
 \end{cases}
\]
Hence, the probability on the left side of Equation~(\ref{eq2.5}) can
be written as
\begin{align}\label{eq2.6}
 &P\{\theta_z \in \Psi_1, \lVert\boldsymbol{z}\rVert < w \} + P\{\theta_z \in \Psi_2, \lVert\boldsymbol{z}\rVert\lvert \mathrm{cos}(\theta_z - \phi/2)\rvert < w \}\notag\\
     &\hphantom{\theta_z \in} + P\{\theta_z \in \Psi_3, \lVert\boldsymbol{z}\rVert\lvert \mathrm{cos}(\theta_z + \phi/2)\rvert < w \}\notag\\
 &=\frac{\phi}{\pi}P\{\lVert\boldsymbol{z}\rVert < w \}\notag\\
     &\hphantom{\theta_z \in}+ 2\int_{\phi/2}^{\pi/2} \frac{1}{2\pi}P\{\lVert\boldsymbol{z}\rVert\lvert \mathrm{cos}(\theta - \phi/2)\rvert < w \}\,\mathrm{d}\theta\notag\\
     &\hphantom{\theta_z \in}+ 2\int_{\pi/2}^{\pi-\phi/2} \frac{1}{2\pi}P\{\lVert\boldsymbol{z}\rVert\lvert \mathrm{cos}(\theta + \phi/2)\rvert < w \}\,\mathrm{d}\theta\notag\\
 &=\frac{\phi}{\pi}P\{R_z < w \}\notag\\
     &\hphantom{\theta_z \in}+ \frac{2}{\pi}\int_{\phi/2}^{\pi/2} P\{R_z\lvert \mathrm{cos}(\theta - \phi/2)\rvert < w \}\,\mathrm{d}\theta\notag\\
  &=\frac{\phi}{\pi} \chi_2^2(w^2)
     + \frac{2}{\pi}\int_{\phi/2}^{\pi/2} \chi_2^2\left(\frac{w^2}{\mathrm{cos}^2(\theta - \phi/2)}\right)\,\mathrm{d}\theta,
 \end{align}
 where $\chi_2^2(.)$ is the cdf of a chi-squared distribution with
 2 degrees of freedom.

\subsection{Method 2: A Method Based on Wynn \& Bloomfield's Approach}
This method calculates $P\{\boldsymbol{z} \in R_2^*\}$ directly, and
is based on the method given by Wynn \& Bloomfield (1971). It can be
seen from Figure 1(b) that the region $R_2^*$ is made up of the circle
with radius $w$ and the remaining area. The probability that
$\boldsymbol{z}$ lies in the circle is given by
\begin{equation}\label{eq2.7}
 P\{\lVert\boldsymbol{z}\rVert < w \} = P\{R_z < w \} = \chi_2^2(w^2)
\end{equation}
The probability that $\boldsymbol{z}$ lies in the remaining region is
4 times the probability that $\boldsymbol{z}$ lies in the slanted-line
shaded area in Figure 1(b), which is
\begin{align}\label{eq2.8}
 &4P\{\theta_z \in [\phi/2,\pi/2],\, \lVert\boldsymbol{z}\rVert > w, \, \lVert\boldsymbol{z}\rVert \mathrm{cos}[\theta_z-\phi/2] < w\}\notag\\
 &=4P\{\theta_z \in [\phi/2,\pi/2],\,w < \lVert\boldsymbol{z}\rVert < w/ \mathrm{cos}[\theta_z-\phi/2] \}\notag\\
 &=4\int_{\phi/2}^{\pi/2}\frac{1}{2\pi}P\{w < \lVert\boldsymbol{z}\rVert < w/ \mathrm{cos}[\theta-\phi/2]\,\mathrm{d}\theta \notag\\
 &=\frac{2}{\pi}\int_{\phi/2}^{\pi/2} \left[\chi_2^2\left(\frac{w^2}{\mathrm{cos}^2[\theta-\phi/2]}\right) - \chi_2^2(w^2) \right]\,\mathrm{d}\theta
\end{align}
Combining Equation~(\ref{eq2.7}) and Equation~(\ref{eq2.8}), we have the following:
\begin{equation}\label{eq2.9}
 P\{\boldsymbol{z} \in R_2^*\} = \frac{\phi}{\pi}\chi_2^2(w^2)
   + \frac{2}{\pi}\int_{\phi/2}^{\pi/2} \left[\chi_2^2\left(\frac{w^2}{\mathrm{cos}^2[\theta-\phi/2]}\right)\right] \,\mathrm{d}\theta,
\end{equation}
which is equal to Equation~(\ref{eq2.6}).
\bigskip

\section{One-Sided Bands}
Similar to the two-sided case, we can express Equation~(\ref{eq2.2})
as
\begin{align}\label{eq3.1}
 &P[\boldsymbol{c'\beta} \leq \boldsymbol{c'\hat{\beta}} + w_u(\boldsymbol{c'F^{-1}c})^{1/2},\,\, \text{for all}\,\, x \in (a,b)]\notag\\
  &= P\left[\sup_{x \in (a,b)}
     \frac{\boldsymbol{(Bc)'z}}{\lVert\boldsymbol{Bc}\rVert} < w_u \right]\notag\\
  &=P\{\boldsymbol{z} \in R_1\},
\end{align}
where $R_1 \subset R^2$ is given by $R_1 = \cap_{x \in (a,b)} R_1(x),$
and
$R_1(x) =
\{\boldsymbol{z}:\boldsymbol{u}'\boldsymbol{z}/\lVert\boldsymbol{u}\rVert
< w\} $
with $\boldsymbol{u} = \boldsymbol{Bc}$.  Rotating the region $R_1$
around the origin the same way as we did in the two-sided case, we
obtain a new region $R_1^*$ with the property
$P\{\boldsymbol{z} \in R_1\} = P\{\boldsymbol{z} \in R_1^*\}$.  The
new region is depicted in Figure 2, and has the expression
\[
  R_1^* = \{\boldsymbol{z}:\boldsymbol{u}'\boldsymbol{z}/\lVert\boldsymbol{u}\rVert < w, \text{for all}\,\,
  \boldsymbol{u} \in E(\phi)\},
\]
where $E(\phi)$ is defined the same as before. Therefore, the
simultaneous confidence level is equal to
\begin{equation}\label{eq3.2}
  P\left[\sup\limits_{\boldsymbol{u} \in E(\phi)}
     \frac{\boldsymbol{u'z}}{\lVert\boldsymbol{u}\rVert} < w_u \right]
 =P\{\boldsymbol{z} \in R_1^*\}.
\end{equation}

\begin{figure}[H]
   \centering
   \caption{Region $R_1^*$}
   \includegraphics[width=4in]{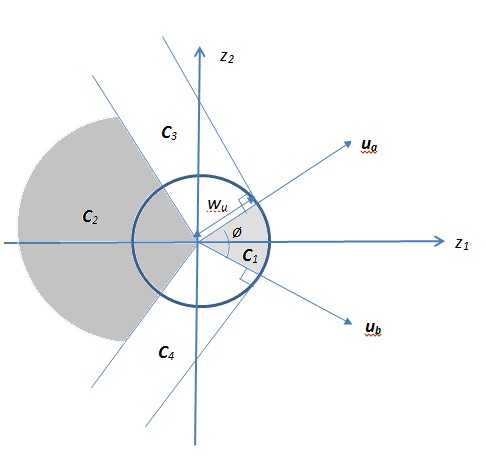}
\end{figure}

\subsection{Method 1: Finding Expression for the Supremum}
Similar to the two-sided bands, this method is based on finding the
supremum in Equation~(\ref{eq3.2}). As before, we write
\[
\frac{\boldsymbol{u'z}}{\lVert\boldsymbol{u}\rVert}
= \lVert\boldsymbol{z}\rVert  \mathrm{cos}(\theta_{uz}),
\]
where $\theta_{uz}$ is the angle between $\boldsymbol{u}$ and $\boldsymbol{z}$. Denote $\Psi_1=[-\phi/2,\phi/2]$,
$\Psi_2=[\phi,-\pi/2]$, and $\Psi_3=[\phi/2,\pi]$. Because of the monotonicity of  $\mathrm{cos}(\theta)$  on $[0, \pi]$, we have
\[
\sup_{\boldsymbol{u} \in E(\phi)} \mathrm{cos}(\theta_{uz}) =
 \begin{cases}
      1                                       &\text{if}\,\, \theta_z \in \Psi_1, \\
      \mathrm{cos}[\theta_z+\phi/2]   & \text{if}\,\, \theta_z \in \Psi_2,\\
      \mathrm{cos}[\theta_z-\phi/2]   & \text{if}\,\, \theta_z \in \Psi_3.\\
 \end{cases}
\]
Hence, the probability on the left side of Equation~(\ref{eq3.2}) can be written as
\begin{align}\label{eq3.3}
 &P\{\theta_z \in \Psi_1, \,\lVert\boldsymbol{z}\rVert < w_u \} + P\{\theta_z \in \Psi_2,\, \lVert\boldsymbol{z}\rVert \mathrm{cos}[\theta_z+\phi/2]  < w_u \}\notag\\
     &\hphantom{\theta_z \in}  +P\{\theta_z \in \Psi_3, \,\lVert\boldsymbol{z}\rVert \mathrm{cos}[\theta_z-\phi/2]  < w_u \}\notag\\
 &=\frac{\phi}{2\pi}P\{\lVert\boldsymbol{z}\rVert < w_u \}\notag\\
     &\hphantom{\theta_z \in}+ 2\int_{\phi/2}^{\pi} \frac{1}{2\pi}P\{\lVert\boldsymbol{z}\rVert \mathrm{cos}[\theta-\phi/2] < w_u\}\,\mathrm{d}\theta \notag\\
 &=\frac{\phi}{2\pi}P\{\lVert\boldsymbol{z}\rVert < w_u \}\notag\\
     &\hphantom{\theta_z \in}+ 2\int_{\phi/2}^{(\pi+\phi)/2} \frac{1}{2\pi}P\{\lVert\boldsymbol{z}\rVert \mathrm{cos}[\theta-\phi/2] < w_u\}\,\mathrm{d}\theta \notag\\
     &\hphantom{\theta_z \in}+ 2\int_{(\pi+\phi)/2}^{\pi} \frac{1}{2\pi}P\{\lVert\boldsymbol{z}\rVert \mathrm{cos}[\theta-\phi/2] < w_u\}\,\mathrm{d}\theta \notag\\
  &=\frac{\phi}{2\pi} \chi_2^2(w_u^2) + \frac{\pi-\phi}{2\pi}
    + \frac{1}{\pi}\int_0^{\pi/2} P\{\lVert\boldsymbol{z}\rVert < \frac{w_u}{\mathrm{cos}(\theta)}\}\,\mathrm{d}\theta.
 \end{align}
 The last step follows from the fact that supremum of
 $\mathrm{cos}(\theta_{uz})$ is negative when
 $\theta_z \in [\pi,(3\pi-\phi)/2] \cup [(\pi+\phi)/2,\pi]$.

\subsection{Method 2: A Method Based on Bohrer \& Francis's Approach}
This method is similar to that for the two-sided case in Section 2.2,
and is based on the method given by Bohrer \& Francis (1972). Notice
that the region $R_1^*$ can be partitioned into four sub-regions,
$C_1, C_2, C_3$, and $C_4$. The probability
$P\{\boldsymbol{z} \in R_1^*\}$ can be calculated by summing up the
probabilities that $\boldsymbol{z}$ lies in the four sub-regions. The
probability that $\boldsymbol{z}$ falls in $C_1$ is equal to
\begin{align}\label{eq3.4}
 P\{\boldsymbol{z} \in C_1\}& = P\{\theta_z \in [(\pi-\phi)/2,(\pi+\phi)/2],\, \lVert\boldsymbol{z}\rVert < w_u\}\notag\\
                            &= \frac{\phi}{2\pi} P\{\lVert\boldsymbol{z}\rVert < w_u\}\notag\\
                            &= \frac{\phi}{2\pi}\chi_2^2(w_u^2).
\end{align}
The probability that $\boldsymbol{z}$ falls into the region $C_2$ is given by
\begin{equation}\label{eq3.5}
P\{\boldsymbol{z} \in C_2\}=P\{\theta_z \in [(\pi+\phi)/2, (3\pi-\phi)/2]\}=\frac{\pi-\phi}{2\pi}.
\end{equation}
Now, if we rotate the region $C_4$ counterclockwise by an angle
$\phi$, then $C_3 \cup C_4$ forms a strip of width $w_u$. By
further rotating the resulting strip clockwise by an angle $\phi/2$, we
obtain a new strip with two sides parallel to the $z_2$ axis. Because
of the rotational invariance of the normal distribution, the
probability that $\boldsymbol{z}$ is in $C_3 \cup C_4$ is equivalent
to the probability of $\boldsymbol{z}$ in the final strip, which is
given by
\begin{equation}\label{eq3.6}
 P\{\boldsymbol{z} \in C_2 \cup C_3\}=P\{0 < Z_1 < w_u\}=\frac{1}{2}\chi_1^2(w_u^2),
\end{equation}
where $\chi_1^2(.)$ is the cdf of a chi-squared distribution with 1
degree of freedom.

Combining Equations~(\ref{eq3.4}), (\ref{eq3.5}), and (\ref{eq3.6}) gives us
\begin{equation}\label{eq3.7}
 P\{\boldsymbol{z} \in R_1^*\}=\frac{\phi}{2\pi}\chi_2^2(w_u^2)+\frac{\pi-\phi}{2\pi}+\frac{1}{2}\chi_1^2(w_u^2).
\end{equation}
It can be shown that Equation~(\ref{eq3.3}) and Equation~(\ref{eq3.7})
are equal.

\section{Monte Carlo Simulation}
Since the confidence bands, defined in Equations~(\ref{eq2.1}) and
(\ref{eq2.2}), are constructed based on the asymptotic properties of
the maximum likelihood estimates for the logistic and probit model
parameters assuming the sample size is large, we conducted Monte Carlo
simulation studies to examine the small-sample performance of the
bands. In our simulation, five different values of
$\boldsymbol{\beta}$ were selected:
$\boldsymbol{\beta}=[-2,\, .3]',\, [0,\, 1.5]',\,[2,\, 5]',\,[-.2,\,
-.3]',\, [-2, -4]'$,
which represent five different forms of the probability response
functions: slowly increasing, moderately increasing, fast increasing,
slowly decreasing, and fast decreasing, respectively.

For the logistic model, we examined three different intervals, narrow,
wide, and extremely wide ("unrestricted"). The endpoints for each
interval were calculated based on the formula
$x=(\text{log}_e[p/(1-p)]-\beta_0)/\beta_1$. For the narrow interval,
$p=.3,\, .7$ were selected for computing the endpoints, and
$p=.1,\, .9$ were selected for computing the endpoints of the wide
interval. We included the very wide interval in the study to compare
the performance of the bands for the restricted and unrestricted
cases, and $p=10^{-10},\, 1-10^{-10}$ were chosen for computing the
endpoints of this interval. The intervals are presented in
Table~\ref{Ta:4.1}.

After the endpoints were calculated for each interval,  we generated values for the
predictor variable $x$ that were equally spread inside each interval with
four different sample sizes ($n = 25, 50, 100, 150$). A uniform
$(0,1)$ random variable was generated, and then the
dichotomous response variable $Y$ was generated based on the
following: $Y=1$ if the uniform random variable was less than
$p(x)=1/(1+\text{exp}(\beta_0+\beta_1x))$; $Y=0$ otherwise.

\begin{table}[H]
 \begin{center}
 \caption{Restricted Intervals For Monte Carlo Simulation}\label{Ta:4.1}
 \begin{tabular}{cccc}
 \hline\noalign{\smallskip}
   $\boldsymbol{\beta}$   &Narrow Interval &Wide Interval  &\textquotedblleft Unrestricted Interval\textquotedblright \\\hline
   $[-2,\, .3]'$   &$(3.842, \, 9.491)$   &$(-.657, \, 13.991)$    &$(-70.086, \, 83.420)$\\
   $[0,\, 1.5]'$   &$(-.565, \, .565)$   &$(-1.465, \, 1.465)$    &$(-15.351, \, 15.351)$\\
   $[2,\, 5]'$   &$(-.569, \, -.231)$   &$(-.839, \, -.039)$    &$(-5.005, \, 4.205)$\\
   $[-.2,\, -.3]'$   &$(-3.491, \, 2.158)$   &$(-7.991, \, 6.657)$    &$(-77.420, \, 76.086)$\\
   $[-2,\, -4]'$   &$(-.712, \, -.288)$   &$(-1.049, \, .049)$    &$(-6.256, \, 5.256)$\\\hline
    \end{tabular}
 \end{center}
 \end{table}

 Three nominal error rates $\alpha$ were considered, $.01, .05$ and
 $.10$.  The computer language \textbf{R}, version 3.2.1, was used to
 carry out all Monte Carlo simulations with $N=5000$ iterations run in
 each simulation to estimate the Monte Carlo error ($1 - $ the
 estimated coverage probability) for the two-sided confidence
 bands. The results are presented in Table~\ref{Ta:4.2}. The results
 in the table suggest that the confidence bands are conservative for
 small samples, and the error approaches the nominal level as the
 sample size increases. Generally, the error reaches the nominal level
 when the sample size is 100, but in some cases it could be as small
 as 50. It is also noted from the table that there is not much
 difference in error between the three different intervals.

\begin{table}[H]
 \begin{center}
 \caption{Estimated Monte Carlo Errors For the Confidence Bands}\label{Ta:4.2}
 \begin{tabular}{c|c|ccc|ccc|ccc}
 \hline\noalign{\smallskip}
   \multirow{2}{*}{$\boldsymbol{\beta}$}
     &\multirow{2}{*}{Sample size $n$}
     &\multicolumn {3}{c|}{Narrow}  &\multicolumn {3}{c|}{Wide}
     &\multicolumn {3}{c}{\textquotedblleft Unrestricted \textquotedblright} \\\cline{3-11}

    & &$\alpha=.01$ &$.05$ &$.10$ &$.01$ &$.05$ &$.10$ &$.01$ &$.05$ &$.10$ \\\hline
  \multirow{4}{*}{$[-2,\, .3]'$}    &25 &$.001$ &$.020$ &$.061$  &$.003$ &$.018$ &$.065$  &$.002$ &$.025$ &$.054$\\
        &50 &$.004$ &$.038$ &$.076$  &$.007$ &$.032$ &$.078$   &$.004$ &$.036$ &$.065$\\
        &100 &$.006$ &$.044$ &$.086$ &$.008$  &$.040$ &$.092$  &$.005$  &$.040$ &$.080$\\
        &150 &$.007$ &$.049$ &$.097$ &$.009$  &$.046$ &$.092$  &$.006$  &$.042$ &$.090$\\\hline
  \multirow{4}{*}{$[0,\, 1.5]'$}    &25 &$.002$ &$.022$ &$.064$  &$.002$ &$.026$ &$.058$  &$.006$ &$.029$ &$.059$\\
        &50 &$.005$ &$.036$ &$.082$  &$.006$ &$.035$ &$.072$  &$.010$ &$.035$ &$.066$\\
        &100 &$.008$ &$.042$ &$.086$ &$.009$ &$.043$ &$.094$  &$.008$ &$.042$ &$.076$\\
        &150 &$.006$ &$.044$ &$.097$ &$.008$ &$.051$ &$.092$  &$.010$ &$.045$ &$.087$\\\hline
  \multirow{4}{*}{$[2,\, 5]'$}    &25 &$.002$ &$.021$ &$.069$ &$.004$ &$.027$ &$.058$  &$.005$ &$.024$ &$.057$\\
        &50 &$.004$ &$.036$ &$.080$  &$.007$ &$.044$ &$.074$   &$.010$ &$.037$ &$.060$\\
        &100 &$.008$ &$.040$ &$.092$ &$.007$ &$.042$ &$.084$   &$.009$ &$.037$ &$.076$\\
        &150 &$.010$ &$.048$ &$.101$  &$.008$ &$.049$ &$.089$  &$.010$ &$.046$ &$.085$\\\hline
  \multirow{4}{*}{$[-.3,\, -.2]'$}    &25 &$.001$ &$.020$ &$.065$  &$.005$ &$.024$ &$.057$  &$.008$ &$.025$ &$.070$\\
        &50 &$.005$ &$.039$ &$.073$  &$.008$ &$.034$ &$.074$  &$.005$ &$.037$ &$.067$\\
        &100 &$.008$ &$.043$ &$.092$ &$.007$ &$.036$ &$.091$  &$.009$ &$.042$ &$.077$\\
        &150 &$.007$ &$.046$ &$.102$ &$.007$ &$.041$ &$.090$  &$.010$ &$.040$ &$.084$\\\hline
  \multirow{4}{*}{$[-4,\, -2]'$}    &25 &$.002$ &$.019$ &$.063$  &$.007$ &$.025$ &$.064$  &$.008$ &$.022$ &$.058$\\
        &50 &$.005$ &$.032$ &$.063$ &$.008$ &$.037$ &$.074$   &$.007$ &$.036$ &$.059$\\
        &100 &$.008$ &$.042$ &$.095$  &$.007$ &$.039$ &$.088$  &$.011$  &$.045$ &$.090$\\
        &150 &$.011$ &$.046$ &$.091$  &$.009$ &$.045$ &$.099$  &$.012$ &$.048$ &$.087$\\\hline
    \end{tabular}
 \end{center}
 \end{table}

 In addition, we considered various other designs.  One case was that
 instead of choosing values that were equally spaced inside the chosen
 intervals for the predictor variable $x$, we chose values for $x$
 that were more concentrated toward one of the endpoints of the
 interval. The values for predictor variable $x$ were generated as
 follows: consider the narrow interval $(3.842, 9.491)$, with $n=25$,
 we first generated $25$ values, stored in a variable $z$, that were
 equally spaced in $(0, 1)$. Then we used the transformation $y = z^6$
 to make the values of $y$ concentrate toward the lower end of the
 interval, 0. Finally, we let $x = 5.649y +3.842$ to obtain the
 desired $x$ values which fell inside the chosen interval and were
 clustered at the lower end, $3.842$.

 The other case examined was when most of the $x$ values were located
 at the center of the chosen intervals, which were generated as
 follows: for the narrow interval $(3.842, 9.491)$ centered at
 $6.6665$, with $n=25$, we first generated $25$ values that were
 equally spaced in $(-1, 1)$, which were stored in a variable called
 $z$. Then we used the transformation $y = z^5$ to make the values of
 $y$ concentrate toward the center of the interval $(-1, 1)$, which is
 0. Finally, we let $x = 5.649/2*y +6.6665$ to obtain the desired $x$
 values which fell inside the chosen interval and were concentrated at
 the center of the interval.

 Since the estimated error rates presented in Table~\ref{Ta:4.2} are
 similar at three different $\alpha$ values, we chose to perform Monte
 Carlo simulations at $\alpha = .05$ for these two additional
 cases. The results are given in Table~\ref{Ta:4.3}, and it is clear
 that the results are consistent with what we observed from
 Table~\ref{Ta:4.2}. We also ran Monte Carlo simulations for the
 one-sided bands in the logistic model, and for both bands in the
 probit model, and again, the results agree with the two-sided case for
 the logistic model; those results are not presented here to avoid
 driving the size of the paper to cumbersome levels.

\begin{table}[H]
 \begin{center}
 \caption{Estimated Monte Carlo Errors For the Confidence Bands at $\alpha=.05$}\label{Ta:4.3}
 \begin{minipage}{\textwidth}
 \begin{tabular}{c|c|cc|cc|cc}
 \hline\noalign{\smallskip}
   \multirow{2}{*}{$\boldsymbol{\beta}$}
     &\multirow{2}{*}{Sample size $n$}
     &\multicolumn {2}{c|}{Narrow}  &\multicolumn {2}{c|}{Wide}
     &\multicolumn {2}{c}{\textquotedblleft Unrestricted \textquotedblright} \\\cline{3-8}

    & &Endpoint \footnote{ The $x$ values are concentrated toward the lower end of the interval} &Center \footnote{ The $x$ values are concentrated toward the center of the interval} &Endpoint &Center &Endpoint &Center\\\hline
  \multirow{4}{*}{$[-2,\, .3]'$}    &25  &$.015$ &$.020$     &$.021$  &$.023$     &$.021$ &$.021$  \\
        &50 &$.022$ &$.031$     &$.031$  &$.038$     &$.036$ &$.040$ \\
        &100 &$.042$ &$.037$     &$.040$  &$.046$     &$.041$ &$.043$ \\
        &150 &$.049$ &$.045$     &$.045$  &$.044$     &$.042$ &$.045$ \\\hline
  \multirow{4}{*}{$[0,\, 1.5]'$}    &25 &$.025$ &$.018$     &$.021$  &$.023$     &$.016$ &$.022$ \\
        &50 &$.032$ &$.031$     &$.039$  &$.038$     &$.025$ &$.030$ \\
        &100 &$.042$ &$.037$     &$.040$  &$.046$     &$.040$ &$.048$ \\
        &150 &$.045$ &$.045$     &$.046$  &$.044$     &$.041$ &$.047$ \\\hline
  \multirow{4}{*}{$[2,\, 5]'$}    &25 &$.022$ &$.018$     &$.015$  &$.022$     &$.016$ &$.021$ \\
        &50 &$.026$ &$.031$     &$.031$  &$.037$     &$.028$ &$.032$ \\
        &100 &$.041$ &$.037$     &$.039$  &$.045$     &$.043$ &$.045$ \\
        &150 &$.045$ &$.045$     &$.042$  &$.043$     &$.042$ &$.048$ \\\hline
  \multirow{4}{*}{$[-.3,\, -.2]'$}    &25 &$.018$ &$.020$     &$.022$  &$.028$     &$.021$ &$.021$ \\
        &50 &$.024$ &$.028$     &$.039$  &$.035$     &$.036$ &$.030$ \\
        &100 &$.041$ &$.043$     &$.038$  &$.042$     &$.045$ &$.044$ \\
        &150 &$.042$ &$.043$     &$.040$  &$.046$     &$.040$ &$.049$ \\\hline
  \multirow{4}{*}{$[-4,\, -2]'$}    &25 &$.018$ &$.020$     &$.016$  &$.028$     &$.021$ &$.025$ \\
        &50 &$.024$ &$.028$     &$.029$  &$.035$     &$.039$ &$.042$ \\
        &100 &$.041$ &$.043$     &$.034$  &$.042$     &$.040$ &$.047$ \\
        &150 &$.041$ &$.044$     &$.042$  &$.046$     &$.042$ &$.049$ \\\hline
    \end{tabular}
  \end{minipage}
 \end{center}
 \end{table}

\section{Example}

If we restrict the independent variable to an interval, the resulting
confidence bands will be narrower than the unconstrained bands.  Here,
we use the data provided by LaVelle (1986) to illustrate the proposed
method. The data are presented in Table~\ref{Ta:5.1}, and these were
the same data that were considered by Piegorsch and Casella
(1988). The study in LaValle (1986) investigated the comutagenic
effects of chromate on frameshift mutagenesis in bacterial assays. The
results presented here are findings for the bacterium
\emph{E}. \emph{coli}, strain 343/435. A control and five doses of the
suspected mutagen, 9-Aminoacridine (9-AA) are reported in
Table~\ref{Ta:5.1}.
\begin{table}[H]
\begin{center}
\captionsetup{justification=centering}
\caption{Mutagenicity of 9-Aminoacridine in \emph{E}. \emph{coli} strains 343/435} \label{Ta:5.1}

\begin{minipage}{11cm}%
\centering
\renewcommand\footnoterule{ \kern -1ex}
\renewcommand{\thempfootnote}{\fnsymbol{mpfootnote}}
\begin{tabular}{c|c|c|c|c|c|c}
\hline
Dose    &-\footnote{The first data pair corresponds to a zero-dose control. The log-dose for this datum was calculated using consecutive-dose average spacing (Margolin \emph{et al.}, 1986).}     &.8   &2.4   &8.0      &24    &80\\\hline
Log-dose      &-1.374    &-.223  &0.875 &2.079   &3.178  &4.382\\\hline
Response       &7/96     &28/96  &64/96 &54/96  &81/96 &96/96\\\hline
\end{tabular}
\vspace{-0.75\skip\footins}
   \renewcommand{\footnoterule}{}
\end{minipage}
\end{center}
\end{table}

As pointed out by Piegorsch and Casella, rather than report confidence
bands over the whole real line, it is often of interest to report
narrower bands over constrained intervals. Furthermore, it is often
noted that human exposure to environmental toxins usually occurs at
low dose levels. Restricting dose levels enables us to direct our
attention to intervals toward the lower end, and hence improve the
confidence limits greatly.

We fit a logistic model to the data using the log-dose level as our
predictor variable $x$, and the ML estimators from the logistic fit
are $\hat{\beta}_0 = -.789$, and $\hat{\beta}_1 = .854$. The inverse of the Fisher
information matrix is
$\boldsymbol{F^{-1}}=\left[ \begin{array}{cc} 0.017 &-0.005 \\ -0.005
    & 0.005 \end{array}\right]$.
If we restrict the predictor variable to an interval $(-1.3,\,.8)$,
then $\boldsymbol{u_a}=(0.163 \quad -0.109)$ ,
$\boldsymbol{u_b}=(0.106 \quad 0.0234)$, and the angle between these
two vectors are $\phi = 0.809$. The two methods for two-sided bands
give the same value for the critical point, $w=2.206$, at a $95\%$
confidence level. Therefore, a $95\%$ two-sided confidence band for
$p(x)$ when $x \in (-1.3, \, .8)$ is given by
\begin{equation*}
 p(x) \in
 \{1 + \text{exp}[(.789-.854x) \pm 2.206(\boldsymbol{c'F^{-1}c})^{1/2}]\} ^{-1}.
\end{equation*}
To make a comparison with the method of Piegorsch and Casella (1988),
the same constrained intervals are considered here, and values of the
critical points $w$ are reported in Table~\ref{Ta:5.2}, along with the
critical points given by Piegorsch and Casella. Values of the critical
point, based on Scheff\'{e}'s method, when there is no restriction on
$x$, are also given here.

\begin{table}[H]
\begin{center}
\captionsetup{justification=centering}
\caption{Critical Values for Two-sided Bands with $\alpha=.05$} \label{Ta:5.2}

\begin{tabular}{c c c}
\hline
Restricted Intervals    & Proposed methods & Piegorsch and Casella  \\\hline
$(-\infty, \,\, \infty)$      &2.447 &2.447\\\hline
(-1.3,\,\,2.0)       &2.344     &2.445 \\\hline
(-1.3,\,\,0.8)       &2.206     &2.274 \\\hline
(-1.3,\,\,-.2)       &2.067     &2.170 \\\hline
\end{tabular}
\end{center}
\end{table}

The sample size in the example is $n=576$, therefore it is sufficiently
large for the asymptotic approximation to hold, based on the Monte
Carlo simulation results in Section~4. It can be seen from the table
that there is clear improvement in the width of confidence bands when
the dose level is restricted to smaller intervals. It is also noted
from the table that the proposed method outperforms Piegorsch and
Casella's procedure by giving smaller critical values and hence
narrower confidence bands. This is what we expected, because as we
pointed out earlier, the embedding procedure used in Piegorsch and
Casella's paper only gave conservative asymptotic confidence bands
when the restricted region for predictor variables is rectangular. In
this data set, the restricted region on the predictor variable, dose
level, is a special case of the rectangular region. The critical
values for one-sided confidence bands over the same intervals are also
presented here in Table~\ref{Ta:5.3}.

\begin{table}[H]
\begin{center}
\captionsetup{justification=centering}
\caption{Critical Values for one-sided bands with $\alpha=.05$} \label{Ta:5.3}

\begin{tabular}{c c c}
\hline
Restricted Intervals    & Proposed methods   \\\hline
(-1.3,\,\,2.0)       &2.049     \\\hline
(-1.3,\,\,0.8)       &1.899     \\\hline
(-1.3,\,\,-.2)       &1.754      \\\hline
\end{tabular}
\end{center}
\end{table}
 Therefore, a $95\%$  upper confidence band for $p(x)$ when $x \in (-1.3, \, .8)$ is given by
\begin{equation*}
 p(x) \leq
 \{1 + \text{exp}[(.789-.854x) - 1.899(\boldsymbol{c'F^{-1}c})^{1/2}]\} ^{-1}.
\end{equation*}

\end{document}